\documentclass[11 pt]{amsart}

\usepackage{amsmath}
\usepackage{amssymb}
\usepackage{amsthm}
\usepackage{caption}
\usepackage{enumitem}
\usepackage{hyperref}
\usepackage{multirow}
\usepackage{nicefrac}
\usepackage{tikz}

\usetikzlibrary{arrows.meta, decorations.markings, decorations.pathmorphing}
\tikzset{->-/.style={decoration={markings,mark=at position .63 with {\arrow{Latex[length=5, width=5]}}}, postaction={decorate}}}
\tikzset{->>-/.style={decoration={markings,mark=at position .45 with{\arrow{Latex[length=5,width=5]}},mark=at position .6 with{\arrow{Latex[length=5,width=5]}}},postaction={decorate}}}

\definecolor{deepblue}{rgb}{0,0,1}
\definecolor{deepred}{rgb}{1,0,0}

\setitemize{noitemsep}

\usetikzlibrary{arrows.meta, decorations.markings, decorations.pathmorphing}
\tikzset{->-/.style={decoration={markings,mark=at position .63 with {\arrow{Latex[length=5, width=5]}}}, postaction={decorate}}}
\tikzset{->>-/.style={decoration={markings,mark=at position .45 with{\arrow{Latex[length=5,width=5]}},mark=at position .6 with{\arrow{Latex[length=5,width=5]}}},postaction={decorate}}}
\tikzset{->>>-/.style={decoration={markings,mark=at position .4 with{\arrow{Latex[length=5,width=5]}},mark=at position .6 with{\arrow{Latex[length=5,width=5]}},mark=at position .8 with{\arrow{Latex[length=5,width=5]}}},postaction={decorate}}}

\definecolor{deepblue}{rgb}{0,0,1}
\definecolor{deepred}{rgb}{1,0,0}
\definecolor{deepgreen}{rgb}{0,1,0}

\newcommand{\Q} {\mathbb{Q}}

\newcommand{\MC} {\operatorname{Mod}}

\newcommand{\Hom} {\operatorname{Hom}}

\newcommand{\Mod}{{\mathcal{M}}}

\newcommand{\Teich}{{\mathcal{T}}}

\newcommand{\Kov}{\mathcal{K}}

\newcommand{\Sym} {\mathbf{S}}

\newcommand{\pt} {\partial}

\newcommand{\wh} {\widehat}

\newcommand{\Lab}{\mathcal{L}}

\newcommand{\Hor}{\mathcal{H}}

\newcommand{\Free}{\mathbf{F}}

\newcommand{\Orb} {\operatorname{Orb}}

\newtheorem{corollary}{Corollary}[section]

\newtheorem{theorem}[corollary]{Theorem}

\newtheorem{lemma}[corollary]{Lemma}

\newtheorem{proposition}[corollary]{Proposition}

\newtheorem{claim}[corollary]{Claim}

\newtheorem{definition}[corollary]{Definition}

\newtheorem{conjecture}[corollary]{Conjecture}

\newtheorem{convention}[corollary]{Convention}

\theoremstyle{remark}

\newtheorem{remark}{Remark}

\begin{document}\title[Random covers]{Tangle free permutations and the Putman-Wieland property of Random covers}

\author[Klukowski and Markovic]{Adam Klukowski and Vladimir Markovi\'c}

\address{\newline Mathematical Institute  \newline University of Oxford  \newline United Kingdom}

\address{\newline All Souls College  \newline University of Oxford  \newline United Kingdom}

\today


\begin{abstract} Let $\Sigma^p_g$ denote a surface of genus $g$ and with $p$ punctures. Our main result is that  the fraction of degree $n$ covers of $\Sigma^p_g$ which have the Putman-Wieland property tends to $1$ as  $n\to \infty$. In addition, we  show that the monodromy  of a random cover of $\Sigma^p_g$  is asymptotically almost surely tangle free.
\end{abstract}

\maketitle

\let\johnny\thefootnote
\renewcommand{\thefootnote}{}

\footnotetext{ This work was supported by the \textsl{Simons Investigator Award}  409745 from the Simons Foundation}
\let\thefootnote\johnny

\section{Introduction}

\subsection{The  Putman-Wieland conjecture}

Let $\Sigma^{p}_{g}$ denote a smooth surface of genus $g\ge 2$ with $p\ge 0$  points removed (which we call cusps or punctures).  Once and for all, we fix a basepoint  $\star \in \Sigma_g^p$ and consider the fundamental group  $\pi_1(\Sigma_g^p,\star)$. 
Denote by  $\MC^{p}_{g}$ the corresponding pure mapping class group. 

By considering  the  basepoint  $\star \in \Sigma_g^p$ as another puncture, we obtain a standard action of $\MC_g^{p+1}$ on $\pi_1(\Sigma_g^p,\star)$.  We associate the following two objects to each  finite index subgroup $K<\pi_1(\Sigma_g^p,\star)$:
\begin{enumerate}
\item Let $\Gamma_K< \MC_g^{p+1}$ denote the  finite index subgroup which leaves $K$ invariant (as a set) with respect to the aforementioned action of $\MC_g^{p+1}$ on $\pi_1(\Sigma_g^p,\star)$. 
\vskip .1cm
\item We say that a pointed cover $\pi:(\Sigma',\star') \to (\Sigma^p_g,\star)$ corresponds to $K$ if $\pi_*\big( \pi_1(\Sigma',\star')  \big)=K$,
where $\pi_*:\pi_1(\Sigma',\star') \to \pi_1(\Sigma^p_g,\star)$ is the induced homomorphism. By $\wh{\Sigma}$ we denote the closed surface obtained by filling in the punctures on $\Sigma'$.
\end{enumerate}

\begin{definition} We say that a finite index subgroup  $K<\pi_1(\Sigma^p_g,\star)$ has the Putman-Wieland property if for each nonzero vector $v \in H^1(\wh{\Sigma},\Q)$, the $\Gamma_K$-orbit  of $v$ is infinite, where $\wh{\Sigma}$ is the  compactification of the corresponding $\Sigma'$.
\end{definition}

Putman and Wieland made the following  conjecture (see Conjecture 1.2. in \cite{p-w}).
\begin{conjecture}\label{conj-pw}
Let $g\ge 2$ and $p\ge 0$. Then every finite index subgroup $K<\pi_1(\Sigma^p_g,\star)$ has  the Putman-Wieland property.
\end{conjecture}

The importance of Conjecture \ref{conj-pw} stems from its close connections with the Ivanov conjectures about the virtual cohomology  of  mapping class groups.  It was shown in \cite{markovic} that Conjecture \ref{conj-pw}  does not hold when $g=2$.

\subsection{A random subgroup has the Putman-Wieland property}

Looijenga \cite{looijenga}, Grunewald-Larsen-Lubotzky-Malestein \cite{g-l-l-m}, Landesman-Litt \cite{l-l-1}, and Markovi\'c-To\v{s}i\'c \cite{m-t},  verified that various types of finite index subgroups of $\pi_1(\Sigma^p_g,\star)$ have the Putman-Wieland property. The purpose of this paper is to prove that  among all  subgroups of index $n$, the fraction of these
which have the Putman-Wieland property  tends to $1$ as $n\to \infty$.

\begin{definition} We let $\Kov_{g,p,n}$  denote the set of  index $n$ subgroups of $\pi_1(\Sigma^p_g,\star)$.  
By $\Kov^{PW}_{g,p,n}$ we denote the subset of $\Kov_{g,p,n}$ consisting of subgroups satisfying the Putman-Wieland property. 
\end{definition}

Our main result says that the fraction of index $n$ subgroups  of $\pi_1(\Sigma^p_g,\star)$ which have the Putman-Wieland property tends to $1$ when $n\to \infty$ provided the genus $g$ is large enough.

\begin{theorem}\label{thm-main-0}
For each $p\ge 0 $ there exists $g_0 \in \mathbb{N}$ such that 
\begin{equation}\label{eq-lim}
\lim\limits_{n\to \infty} \frac{|\Kov^{PW}_{g,p,n}|}{|\Kov_{g,p,n}|}=1,
\end{equation}
when $g\ge g_0$. 
\end{theorem}
In fact, we prove a  stronger statement: 

\begin{theorem}\label{thm-main} For each $\kappa < \frac{1}{2}$ there exists $g_0 \in \mathbb{N}$ such that  (\ref{eq-lim}) holds when $g\ge g_0$, and $p \leq g^{\kappa}$. 
\end{theorem}

\subsection{Random permutations are tangle free}

The Symmetric group $\Sym_n$ is the group of permutations of the set $[n]=\{1,\dots,n\}$. Let $\Free_m$ denote the free group on $m$ generators, and  $\Hom_{m,n}$ the set of homomorphisms from $\Free_m$ to the symmetric group $\Sym_n$.

\begin{definition}\label{definition-tang}  Let $R>0$. We say that $w_1,w_2\in \pi_1(\Sigma^p_g,\star)$ are $R$-tangled by $\phi\in \Hom_{m,n}$  if there exists $k\in [n]$ such that 
$$
\left|\Orb(\phi(w_1),k)\right| +  \left|\Orb(\phi(w_2),k)\right|\le R,
$$
where $\Orb(\phi(w_i),k)\subset [n]$ is the orbit of $k$ under the iterates of $\phi(w_i)$.
\end{definition}

One of the key ingredients in the proof of Theorem \ref{thm-main} is the following theorem.

\begin{theorem}\label{thm-comb-0}
Let $w_1, w_2 \in \Free_m$ be two elements  whose nontrivial powers are all distinct. Then for every $R>0$ the equality
$$
\frac{\left|\{ \phi\in \Hom_{m,n} : w_1,w_2 \,\, \text{are R-tangled by} \,\, \phi \}\right|}{\left|\Hom_{m,n}\right|} \leq \frac{C}{n}
$$
holds, where $C$ depends only on $w_1, w_2, R$, and $m$.
\end{theorem}
\begin{remark}
It was shown by Monk-Thomas \cite{m-t} that for any $R>0$ there exists $g_0$ such that a random Riemann surface 
of genus $g>g_0$, picked with respect to the Weil-Petersson probability measure, is $R$-tangle-free.  Applying Theorem \ref{thm-comb-0} one can prove a version of this result for random covers of cusped hyperbolic surfaces.
\end{remark}

\subsection{Organisation of the paper} Each section starts with a brief outline. Here we only give a  broad overview of the paper. 
In Section \ref{section-druga} we state the result by Markovi\'c-To\v{s}i\'c \cite{m-t} that covers admitting sufficiently large spectral gap  have the Putman-Wieland property. This naturally leads us to the results by Magee-Naud-Puder \cite{m-n-p}, and Hide-Magee \cite{h-m}, that a random cover of a fixed hyperbolic surface has no new (sufficiently) small eigenvalues. The combination of these two results is the main idea behind the proof of Theorem \ref{thm-main}. 

The main difficulty we need to overcome is be able to compare the spectral gaps of a random cover and its compactification (obtained by filling in the punctures). We do this by showing that a random cover of a cusped hyperbolic surface has the $L$-horoball property, which by the work of Brooks \cite{brooks} guarantees that the two spectral gaps are comparable. Assuming Theorem \ref{thm-adam} we prove Theorem \ref{thm-main} in Section \ref{section-horo}.

The remainder of the paper (after Section \ref{section-horo}) is mostly devoted to proving Theorem \ref{thm-adam}.
As explained in Section \ref{section-L},  showing that a random cover has the $L$-horoball property reduces to proving Theorem \ref{thm-comb-0}.  The proof of Theorem \ref{thm-comb-0} has combinatorial flavour and we explain the main steps in Sections \ref{section-comb} and \ref{section-comb-1} where we present the proof.

\section{Spectral gap and the Putman-Wieland property}\label{section-druga}
In this section we recall results  which underpin the proof of Theorem \ref{thm-main}.
The first  is that a subgroup $K\in\Kov_{g,p,n}$ has the Putman-Wieland property assuming that the spectral gap of the surface $\wh{X}_K$ is uniformly bounded away from zero. Here $\wh{X}_K$ denotes the compactification  of the corresponding   holomorphic covering surface  $X_K$.

Next, we define the subset $\Kov^{\lambda_{1}}_{n}(X) \subset \Kov_{g,p,n}$ consisting of  subgroups $K\in \Kov_{g,p,n}$ for which $X_K$ has the same spectral gap as fixed hyperbolic surface $X$. We then state the results from \cite{m-n-p}, \cite{h-m} saying that a random element of  $\Kov_{g,p,n}$ is asymptotically almost surely in $\Kov^{\lambda_{1}}_{n}(X)$ (for suitable $X$).

\subsection{Pointed holomorphic covers}\label{subsection-pointed}  Let $\Mod_{g,p}$  denote the moduli space of Riemann surfaces of genus $g$ and with $p$ cusps. A marked pointed Riemann surface is a triple $(X,x,f)$, where $X\in \Mod_{g,p}$, $x\in X$, and $f:(\Sigma^p_g,\star)\to (X,x)$ is pointed homeomorphism. 

Let $X\in \Mod_{g,p}$, and choose  a marked pointed Riemann surface  $(X,x,f)$.
Fix $K \in \Kov_{g,p,n}$, and denote by $\pi:(\Sigma',\star') \to (\Sigma^p_g,\star)$ a pointed cover corresponding to $K$. Then there exist a unique   marked pointed Riemann surface $(X_K,x_K,f_K)$,  and  a holomorphic unbranched covering  $\pi_K:X_K \to X$, such that the following diagram commutes
\begin{center}
\begin{tikzpicture}[every node/.style={midway}]
  \matrix[column sep={7em,between origins}, row sep={2em}] at (0,0) {
    \node(R) {$(X_K,x_K)$}  ; & \node(S) {$(\Sigma',\star')$}; \\
    \node(R/I) {$(X,x)$}; & \node (T) {$(\Sigma^p_g,\star)$};\\
  };
  \draw[<-] (R/I) -- (R) node[anchor=east]  {$\pi_K$};
  \draw[<-] (R) -- (S) node[anchor=south] {$f_K$};
  \draw[->] (S) -- (T) node[anchor=west] {$\pi$};
  \draw[<-] (R/I) -- (T) node[anchor=north] {$f$};
\end{tikzpicture}
\end{center}
Note that the Riemann surface $X_K$ depends only on $K$ and not on  the choice of the point $x\in X$ or the marking $f$. 
\begin{remark}
On the other hand, the point $x_K$ and the map $f_K$ depend  on both $x$ and $f$.
\end{remark}
We close this subsection with the following definition.
\begin{definition} 
By  $\wh{X}_K$ we denote the closed Riemann surface obtained from $X_K$ by filling in the punctures.
\end{definition}

\subsection{A random cover retains the spectral gap}

The starting point in the proof of Theorem \ref{thm-main} is  the following result by 
Markovi\'c-To\v{s}i\'c  \cite{m-t} which states that covers admitting sufficiently large spectral gap must have the Putman-Wieland property (see Theorem 1.9 in \cite{m-t}).

\begin{theorem}\label{thm-mt} Let $g\ge 2$, $p\ge 0$, and  $K \in \Kov_{g,p,n}$. Suppose that there exists  $X\in \Mod_{g,p}$ such that 
$$
\frac{1+2\lambda_1(\wh{X}_K)}{2\lambda_1(\wh{X}_K)}\le g.
$$
Then $K \in \Kov^{PW}_{g,p,n}$. Here $\lambda_1$ denotes the smallest non-zero eigenvalue of the (hyperbolic) laplacian.
\end{theorem}

This result indicates that in order to show $K \in \Kov^{PW}_{g,p,n}$ it suffices to bound from below $\lambda_1(\wh{X}_K)$ for some $X\in \Mod_{g,p}$ (this bound needs to be uniform in $n$). This brings us to the second key ingredient in the proof Theorem \ref{thm-main} which is the result that a random cover of $X$ does not have any new small eigenvalues.

\begin{definition} For $X\in \Mod_{g,p}$, we set
$$
\Kov^{\lambda_{1}}_{n}(X)= \big\{ K \in \Kov_{g,p,n} :\,\, \lambda_1 (X_K) = \lambda_1(X) \big\}.
$$
\end{definition}
Note that the inequality   $\lambda_1 (X_K) \le  \lambda_1(X)$ always holds because the pullback of an
eigenfunction on $X$  is  an eigenfunction on $X_K$. Requiring  $\lambda_1 (X_K) =  \lambda_1(X)$ means that the laplacian on the covering surface $X_K$ has no new eigenvalues which are strictly smaller than $\lambda_1(X)$. The second key ingredient in the proof of Theorem \ref{thm-main} is the following:

\begin{theorem}\label{thm-magee}  Let $X \in \Mod_{g,p}$ be such that $\lambda_1(X)< \frac{3}{16}$. Then 
\begin{equation}\label{eq-magee} 
\lim\limits_{n\to \infty} \frac{|\Kov^{\lambda_{1}}_{n}(X)|}{|\Kov_{g,p,n}|}=1.
\end{equation}
\end{theorem}
This theorem follows from deep results by Magee-Naud-Puder \cite{m-n-p} in the case $p=0$, and by Hide-Magee \cite{h-m} when $p>0$. However, they prove this with respect to the uniform measure on  the degree $n$ covers of $X$ with a labelled fibre (which is denoted by $\Lab_{g,p,n}$ in Section \ref{section-labelled}), while Theorem \ref{thm-magee} is  the version of their results with respect to the uniform measure on the set $\Kov_{g,p,n}$ of index $n$ subgroups of $\pi_1(\Sigma^p_g,\star)$. The two models are closely related which we explain in Section \ref{section-labelled} (where we formally prove Theorem  \ref{thm-magee}).  

It is clear that combining Theorem \ref{thm-mt} and Theorem \ref{thm-magee} brings us closer to proving Theorem \ref{thm-main}. However, the main  obstacle is that the statement of  Theorem \ref{thm-mt} inputs 
$\lambda_1(\wh{X}_K)$, while the statement of Theorem \ref{thm-magee} outputs $\lambda_1(X_K)$. Therefore, we have to show that for a random cover of $X$ these  two geometric quantities are in some sense related to each other.

\section{The $L$-horoball property}\label{section-horo}

We explain the notion of the $L$-horoball property and  its significance in relating the spectral gaps of a cusped hyperbolic surface and its compactification. Then we define the subset $ \Kov^{\mathcal{H}}_n(X,L)\subset \Kov_{g,p,n}$ consisting of  subgroups for which the induced holomorphic covering surface  $X_K$ has the $L$-horoball property. At the end of the section we prove 
Theorem \ref{thm-main} assuming Theorem \ref{thm-adam} which states that a random subgroup belongs to  $\Kov^{\mathcal{H}}_n(X,L)$ asymptotically almost surely. 

\begin{remark}
In \cite{b-m}, Brooks and Makover  developed a certain model of random surfaces and  proved a random surface in this model has the $L$-horoball property for every $L$. 
\end{remark}

\subsection{The $L$-horoball property and the Cheeger constant}
In this subsection we let $S$ denote a hyperbolic surface. To control the behaviour of $\lambda_1$ under conformal compactification we use the $L$-horoball property devised by Brooks \cite{brooks}. 

\begin{definition} Let  $S$ be a hyperbolic Riemann surface with at least one cusp.  Given $L>0$, we say that $S$ has the $L$-horoball property if the  horoballs  of perimeter $L$ around all punctures are pairwise disjoint and embedded.
\end{definition}
\begin{remark}
It is a known feature of hyperbolic geometry that every hyperbolic cusped surface has the $1$-horoball property.
\end{remark}
Let $S^c$ denote the closed surface obtained from $S$ by filling in the puncture. Brooks (see Theorem 4.1 in \cite{brooks}) established a connection between the Cheeger constants of $S$ and $S^c$ which we denote by $h(S)$ and $h(S^c)$  respectively.
\begin{theorem}\label{thm-brooks} For every $C>1$ there exists $L>0$ such that if $S$ is a finite area hyperbolic surface which has the $L$-horoball property then
\begin{equation}\label{eq-brooks}
\frac{1}{C}h(S)\le h(S^c) \le Ch(S). 
\end{equation}
\end{theorem}

\subsection{The $L$-horoball property and the spectral gap}

Let us define the set of covers with the $L$-horoball property.

\begin{definition} For $L>0$, and $X\in \Mod_{g,p}$, we let
$$
\Kov^{\mathcal{H}}_{n}(X,L)= \big\{ K \in \Kov_{g,p,n} : \,\, X_K \,\,\, \text{has L-horoball property}\big\}.
$$
\end{definition}

We  prove the following lemma by combining   Theorem \ref{thm-brooks} with the classical results by Cheeger and Buser.

\begin{lemma}\label{lemma-sve} There exist universal constants $q,L>0$ such that  for every 
$X \in \Mod_{g,p}$, and every   $K \, \in \, \Kov^{\lambda_{1}}_{n}(X)\cap \Kov^{\mathcal{H}}_n(X,L)$, the inequality
\begin{equation}\label{eq-sve-0} 
\lambda_1 (\wh{X}_K) \ge q \lambda^2_1(X)
\end{equation}
holds. 
\end{lemma}

\begin{proof}  Let $L$ be the constant from Theorem \ref{thm-brooks} such that the inequalities (\ref{eq-brooks}) hold for $C=2$. Consider any $K \, \in \, \Kov^{\lambda_{1}}_{n}(X)\cap \Kov^{\mathcal{H}}_n(X,L)$. 
We show that there exists a universal constant $q>0$ such that (\ref{eq-sve-0}) holds.

Firstly, the classical Cheeger's inequality gives 
\begin{equation}\label{eq-jutro}
\lambda_1(\wh{X}_K) \geq \frac{1}{4} h^2(\wh{X}_K).
\end{equation}
Moreover, from the choice of  $L$, and the assumption   
$K \in \Kov^{\mathcal{H}}_n(X,L)$, we get 
\begin{equation}\label{eq-jutro-1}
h(\wh{X}_K) \geq \frac{1}{2} h(X_K).
\end{equation}
Furthermore, by Theorem 7.1 of \cite{buser-0} there exists an absolute constant $c>0$ such that 
\begin{equation}\label{eq-jutro-1}
h (X_K) > c \lambda_1 (X_K).
\end{equation}
Putting together the last three inequalities shows that 
$$
\lambda_1(\wh{X}_K)\ge \frac{c^2}{16} \lambda^2_1 (X_K).
$$
This inequality, combined with the assumption  $K \in \Kov^{\lambda_{1}}_{n}(X)$, yields
$$
\lambda_1(\wh{X}_K)\ge \frac{c^2}{16} \lambda^2_1 (X).
$$
This proves (\ref{eq-sve-0}) for  $q=\tfrac{c^2}{16}$. 
\end{proof}

\subsection{A random cover has the $L$-horoball property}

The following theorem shows that a random cover of $X$ has the $L$-horoball property. It is proved in Section \ref{section-L}.

\begin{theorem}\label{thm-adam}  For every $X \in \Mod_{g,p}$, and every $L>0$, we have
\begin{equation}\label{eq-adam} 
\lim\limits_{n\to \infty} \frac{|\Kov^{\Hor}_{n}(X,L)|}{|\Kov_{g,p,n}|}=1.
\end{equation}
\end{theorem}

The equality (\ref{eq-adam}) puts us in the position to put together  Theorem \ref{thm-mt} and Theorem \ref{thm-magee} to 
prove the main result Theorem \ref{thm-main}.

\subsection{Proof of Theorem \ref{thm-main}}
We begin by stating the  result which follows from the paper by Hide \cite{hide}.
\begin{proposition}\label{prop-hide}
For every $\kappa < \frac{1}{2}$ there exist constants $g_1> 0 $, and $0<\delta <\tfrac{3}{16}$, with the following property. If $g \geq g_1$ and $p \leq g^{\kappa}$, then there exists  $X \in \Teich(\Sigma_{g, p})$ such that 
\begin{equation}\label{eq-hide}
\delta\le \lambda_1 (X) <\frac{3}{16}.
\end{equation}
\end{proposition}
\vskip .3cm

Fix $\kappa < \frac{1}{2}$, and let $q$ and $L$ be the constant  from Lemma \ref{lemma-sve}. Set
$$
g_0=\max\left\{g_1, \frac{1+2q\delta^2}{2q\delta^2} \right\}.
$$

\begin{claim}\label{claim-marko} Suppose $g\ge g_0$,  $p\le g^\kappa$. Then there exists $X\in \Mod_{g,p}$ such that 
$\Kov^{\lambda_{1}}_{n}(X)\cap \Kov^{\mathcal{H}}_n(X,L) \, \subset \, \Kov^{PW}_{g,p,n}$.
\end{claim}

\begin{proof} Let $X\in \Mod_{g,p}$ be such that (\ref{eq-hide}) holds, and suppose 
$K\, \in \, \Kov^{\lambda_{1}}_{n}(X)\cap \Kov^{\mathcal{H}}_n(X,L)$. Then by the inequality (\ref{eq-sve-0}), and from (\ref{eq-hide}) , we derive the inequality
$$
\lambda_1 (\wh{X}_K) \ge q \delta^2. 
$$
Combining this with the lower bound 
$$
g\ge \frac{1+2q\delta^2}{2q\delta^2}, 
$$
and applying Theorem \ref{thm-mt}, proves that $K \in \Kov^{PW}_{g,p,n}$. 
\end{proof}

To finish the proof, we first observe the equality 
\begin{equation}\label{eq-sve} 
\lim\limits_{n\to \infty} \frac{|\Kov^{\lambda_{1}}_{n}(X)\cap \Kov^{\mathcal{H}}_n(X,L) |}{|\Kov_{g,p,n}|}=1.
\end{equation}
This follows by from Theorem \ref{thm-adam} and Theorem \ref{thm-magee}.  The proof of Theorem \ref{thm-main} now follows from  the equality (\ref{eq-sve}) and Claim \ref{claim-marko}.

\section{Finite covers with labelled fibres}\label{section-labelled}

In this section we define the set $\Lab_{g,p,n}$ of fibre labelled covers of $\Sigma^p_g$ which consists of monodromy homomorphisms
of (fibre unlabelled) covers  $\Kov_{g,p,n}$. We then observe that the natural projection $P_n:\Lab_{g,p,n} \to \Kov_{g,p,n}$ enables us to easily replace $\Kov_{g,p,n}$ by $\Lab_{g,p,n}$ in the statements of Theorem \ref{thm-magee} and Theorem \ref{thm-adam}. 
At the end of the section we derive the proof of Theorem \ref{thm-magee}, and state  Theorem \ref{thm-comb} which is a version of Theorem \ref{thm-comb-0} for transitive homomorphisms.

\subsection{Labelled fibres and the Symmetric group}

We say that a pair $(\pi,\iota)$ is a degree $n$ cover with a labelled fibre if
\begin{enumerate}
\item $\pi:\Sigma'\to \Sigma^p_g$  is a (connected)  cover of degree $n$,
\vskip .1cm
\item $\iota:\pi^{-1}(\star)\to [n]$ a labelling
\end{enumerate}
(recall the abbreviation $[n]=\{1,\dots,n\}$). Two such covers $\pi':\Sigma'\to \Sigma^p_g$, and  $\pi'':\Sigma''\to \Sigma^p_g$, are equivalent if there exists a homeomorphism $I:\big(\Sigma',(\pi')^{-1}(\star)\big)\to \big(\Sigma'',(\pi'')^{-1}(\star)\big)$ so that $\pi'=\pi''\circ I$, and $I \circ (\iota'')^{-1}= (\iota')^{-1}$.

\begin{definition}
The set of equivalence classes of degree $n$ covers with a labelled fibre is denoted by $\Lab_{g,p,n}$.
\end{definition}
To each equivalence class $[\pi,\iota]\in \Lab_{g,p,n}$ we associate the monodromy homomorphism $\phi_{[\pi,\iota]}:\pi_1(\Sigma^p_g, \star)\to \Sym_n$ which describes how the fibre $\pi^{-1}(\star)$ is permuted when following lifts of a closed loop from $\Sigma^p_g$ to $\Sigma'$. 
In fact, the equivalence class  $[\pi,\iota]$  is uniquely determined by the  monodromy homomorphism $\phi_{[\pi,\iota]}$ (see Section 1 in \cite{m-n-p}).

Since the cover $\pi$ is connected it follows that the homomorphism $\phi_{[\pi,\iota]}$ is transitive (i.e. the image group 
$\phi_{[\pi,\iota]}\big(\pi_1(\Sigma^p_g, \star)\big)$ acts transitively on the set $[n]$).
 Therefore, there is a natural bijection
\begin{equation}\label{eq-bijection}
\Lab_{g,p,n}\,\, \longleftrightarrow \,\, \{\text{transitive homomorphisms $\pi_1 (\Sigma^p_g, \star) \rightarrow \Sym_n$}  \}.
\end{equation}
\vskip .5cm
\begin{convention}\label{convention}
Using  bijection (\ref{eq-bijection}), we let $\Lab_{g,p,n}\subset \Hom_{m,n}$  denote the set of transitive homomorphisms $\pi_1 (\Sigma^p_g, \star) \to \Sym_n$. 
\end{convention}
\vskip .5cm

\begin{remark} It follows from the work by Liebeck-Shalev \cite{l-s} (which generalises an old theorem of Dixon \cite{dixon}, see also  the introduction in \cite{p-z}) that non-transitive homomorphisms $\pi_1 (\Sigma^p_g, \star) \rightarrow \Sym_n$ are statistically insignificant when $n$ is large for any fixed $g$ and $p$ such that $3g+p-3>0$. We have
\begin{equation}\label{eq-limit}
\lim_{n\to\infty} \frac {\left| \Lab_{g,p,n}\right|}{\left| \Hom_{m,n} \right|}= 1. \end{equation}
\end{remark}

\subsection{Subgroups and pointed covers} Two pointed covers $\pi':(\Sigma',\star')\to (\Sigma^p_g,\star)$, and $\pi'':(\Sigma'',\star'')\to (\Sigma^p_g,\star)$, are equivalent if there exists a pointed homeomorphism $\text{I}:(\Sigma',\star')\to (\Sigma'',\star'')$ such that $\pi'=\pi''\circ I$. The equivalence class  of a pointed cover $\pi:(\Sigma',\star') \to (\Sigma^p_g,\star)$ is uniquely determined by the subgroup $\pi_*\big( \pi_1(\Sigma',\star')  \big)$, where $\pi_*:\pi_1(\Sigma',\star') \to \pi_1(\Sigma^p_g,\star)$ is the induced homomorphism.
Therefore, the set  $\Kov_{g,p,n}$ is in the bijection with the set of equivalence classes of degree $n$ pointed covers of 
$\Sigma^p_g$.

Let $(\pi,\iota)$ be a fibre labelled cover $\pi:\Sigma'\to \Sigma^p_g$. Set $\star'=\iota^{-1}(1)$, and consider the pointed cover
$\pi:(\Sigma',\star')\to (\Sigma^p_g,\star)$. It is elementary to check that to equivalent fibre labelled covers we associate equivalent pointed covers. Thus, we have constructed the map
\begin{equation}\label{eq-lab}
P_n:\Lab_{g,p,n} \to \Kov_{g,p,n}.
\end{equation}
Moreover, if  $\phi,\psi\in \Lab_{g,p,n}$ then $P_n(\psi)=P_n(\phi)$ if and only if the homomorphism $\phi$ and $\psi$ agree 
up to post-conjugation by a permutation  of the set $\{2,\dots,n\} $.
There are exactly $(n-1)!$ such permutations. This enables us to conclude:

\begin{lemma}\label{lemma-P} The pre-image (under the map $P_n$) of each element of  $\Kov_{g,p,n}$ consists of exactly $(n-1)!$ different elements of $\Lab_{g,p,n}$.
\end{lemma}
Thus, each index $n$ subgroup $K\in \Kov_{g,p,n}$ corresponds to exactly $(n-1)!$  homomorphisms from $\Lab_{g,p,n}$. 

\begin{definition} Let $\phi\in\Lab_{g,p,n}$. For $X\in \Mod_{g,p}$, we let $X_\phi=X_K$, where $K=P_n(\phi)$. 
\end{definition} 

It is important to observe that $X_\phi=X_\psi$ if $P_n(\phi)=P_n(\psi)$.
Explicitly, the subgroup $K$ is the stabiliser of $1$ in the action $\phi\big(\pi_1(\Sigma^p_g,\star)\big)$  on $[n]$.

\subsection{Proof of Theorem \ref{thm-magee} } Suppose  $X\in \Mod_{g,p}$ is such that $\lambda_1(X)<\tfrac{3}{16}$. 
Let $\Lab^{\lambda_{1}}_{n}(X)\subset \Lab_{g,p,n}$ denote the set of homomorphisms $\phi$ such that $\lambda_1(X)=\lambda_1(X_\phi)$.
It was shown in \cite{m-n-p} for $p=0$, and in \cite{h-m} when $p>0$, that 
\begin{equation}\label{eq-lamb}
\lim\limits_{n\to \infty} \, \frac{\left| \Lab^{\lambda_{1}}_n(X)  \right|}{\left| \Lab_{g,p,n}  \right|}=1.
\end{equation}
On the other hand, it follows from Lemma \ref{lemma-P}  that
$$
\left|\Lab^{\lambda_{1}}_{n}(X)\right|= (n-1)!\left|\Kov^{\lambda_{1}}_{n}(X)\right|,\,\,\,\,\,\,\,\,\,\, \left|\Lab_{g,p,n}\right|= (n-1)!\left|\Kov_{g,p,n}\right|.
$$
Together with (\ref{eq-lamb}) this  yields the proof of Theorem \ref{thm-magee}.

\subsection{Random transitive permutations are tangle free}\label{subsection-tang}

The proof of Theorem \ref{thm-adam} consists of a geometric and a combinatorial part. The combinatorial part reduces to the statement that a transitive random permutation is asymptotically almost surely tangle free. 

\begin{theorem}\label{thm-comb}
Suppose $p>0$, and let $w_1, w_2 \in \pi_1(\Sigma^p_g)$ denote two elements  whose nontrivial powers are all distinct. Then for every $R>0$, the equality
$$
\frac{\left|\{ \phi\in \Lab_{g,p,n} : w_1,w_2 \,\, \text{are R-tangled by} \,\, \phi \}\right|}{\left|\Lab_{g,p,n}\right|} \leq \frac{C}{n}
$$
holds, where $C$ depends only on $w_1, w_2, R, g$, and $p$.
\end{theorem}
\begin{proof} 
The combination of the equality (\ref{eq-limit}) and Theorem \ref{thm-comb-0} yields the proof of Theorem \ref{thm-comb}. 
\end{proof}

\section{Generic  covers have $L$-horoball property}\label{section-L}
In this section we explain the geometric content of the proof of Theorem \ref{thm-adam} and conclude its proof assuming Theorem \ref{thm-comb}.

\begin{definition}
Given $X\in \Mod_{g,p}$, we let $\Lab^{\Hor}_{n}(X,L) \subset \Lab_{g,p,n}$ denote the set of homomorphisms $\phi$ such that $X_\phi$ has the $L$-horoball property. 
\end{definition}
In view of Lemma \ref{lemma-P}, to prove 
Theorem \ref{thm-adam} it suffices to prove the following:
\begin{theorem}\label{thm-adam-1}  For every $X \in \Mod_{g,p}$, and every $L>0$, we have
\begin{equation}\label{eq-adam-1} 
\lim\limits_{n\to \infty} \frac{|\Lab^{\Hor}_{n}(X,L)|}{|\Lab_{g,p,n}|}=1.
\end{equation}
\end{theorem}

The remainder of the paper is devoted to proving Theorem \ref{thm-adam-1}. In this section we reformulate the $L$-horoball property in terms of the pairs of connected cusps on $X$ which behave well under covers. This allows us to prove Theorem \ref{thm-adam-1}  assuming Theorem \ref{thm-comb}.

Below we prove  Propositions \ref{prop-pair-1} and \ref{prop-pair-2} which relate the $L$-horoball property of a covering $\pi_\phi:X_\phi\to X$ to the branching degrees at pairs of connected cusps. Then, in Proposition \ref{prop-pair-5} we  express these branching degrees  in terms of the combinatorial data of the homomorphism $\phi\in \Lab_{g,p,n}$.

\subsection{Geometry of pairs of connected cusps}
In the remainder of this section we assume $X\in \Mod_{g,p}$ and $p>0$. 
\begin{definition}
A pair of connected cusps on  $X$ is a triple $(c_1,c_2 ,\delta)$, where $c_1, c_2$ are cusps on $X$, and $\delta$ is a geodesic joining them.
\end{definition}

Note that we allow $c_1=c_2$. We will need some more vocabulary to work with pairs of connected cusps. If $c$ is a cusp on $X$, 
we let $\Hor_c(r)$ denote the horoball  based at $c$ such that that the horocycle $\pt{\Hor_c(r)}$ has perimeter equal to $r$.
\begin{definition} Let $(c_1, c_2, \delta)$ be a pair of connected cusps on $X$.
\begin{itemize}
\item The beam  $\beta (c_1,c_2,\delta)$ is  the segment of $\delta$ lying outside the  horoballs $\Hor_{c_{1}} (1)$ and $\Hor_{c_{2}} (1)$. 
\vskip .1cm
\item  The length of $(c_1, c_2, \delta)$ is denoted by  $|(c_1, c_2, \delta)|$. It is defined as the  length of the beam $\beta (c_1, c_2, \delta)$.\end{itemize}
\end{definition}

Having defined pairs of connected cusps we  now state a few of their  properties. We begin with the following elementary fact from hyperbolic geometry.
\begin{lemma}\label{lemma-pomoc} Let $c$ be a cusp of $X$. Suppose $r\ge 1$. Then the distance  between the horocycles
$\pt{\Hor_c(r)}$ and  $\pt{\Hor_c(1)}$ is $\log r$. 
\end{lemma}
The following proposition explains  the relationship between pairs of connected cusps and  the $L$-horoball property. 
\begin{proposition}\label{prop-pair-1}
A hyperbolic surface $X$ has the $L$-horoball property if and only if the length of every pair of connected cusps on $X$ is at least $2 \log L$.
\end{proposition}
\begin{proof} Suppose $c_1$ and $c_2$ are cusps on $X$ (not necessarily distinct). From Lemma \ref{lemma-pomoc} we conclude
that the horoballs  $\Hor_{c_{1}} (L)$ and $\Hor_{c_{2}} (L)$ are disjoint  if and only if for each geodesic $\delta$ connecting $c_1$ and $c_2$, the length $|(c_1, c_2, \delta)|$ is at least $2 \log L$. This proves the proposition.
\end{proof}

The reason why it is convenient to reformulate the $L$-horoball property in terms of pairs of connected cusps is that they behave well under covering maps.  

\begin{proposition}\label{prop-pair-2} Suppose  $\pi:X'\to X$ is a holomorphic covering. Let  $(c'_1, c'_2, \delta')$ be a pair of connected cusps on $X'$, and $(c_1, c_2, \delta)$ a pair of connected cusps on $X$, such that $\pi(c'_1, c'_2, \delta')=(c_1, c_2, \delta)$. Denote the branching degrees of $\pi$ at $c'_1, c'_2$, by $d_1$ and  $d_2$ respectively. 
Then 
\begin{equation}\label{eq-pair-2}
|(c'_1, c'_2, \delta')| = |(c_1, c_2, \delta)| + \log d_1 + \log d_2.
\end{equation}
\end{proposition}
\begin{proof}
The preimage of $\Hor_{c_i}(1)$ is the horoball $\Hor_{c'_i} (d_i)$ on which $\pi$ is a $d_i$-fold cyclic covering. From Lemma \ref{lemma-pomoc} we find that the length of the segment $\beta(c'_1,c'_2,\delta')$ is given by
$$
|\beta(c'_1,c'_2,\delta')|=|\beta(c_1,c_2,\delta)|+\log d_1 + \log d_2,
$$
which proves the proposition.
\end{proof}

Finally, let us observe that the set of pairs of connected cusps of bounded  length is finite.

\begin{proposition}\label{prop-pair-3}
Let $R>0$, and fix $X\in \Mod_{g,p}$. Then the set of pairs of connected cusps whose length is at most $R$ is finite.
\end{proposition}
\begin{remark} It can be shown that the number of such pairs of connected cusps is at most $2e^R\times ( \text{the number of cusps on X})$.
\end{remark}
\begin{proof}  The set of pairs of connected cusps on $X$ whose length is at most $R$ is both compact and discrete, and thus finite. 
\end{proof}

\subsection{Lollipops} 
To each pair of connected cusps $(c_1, c_2, \delta)$ on $X$ we associate two elements of the fundamental group of $X$. 
By $x=x (c_1, c_2, \delta)$ we denote the midpoint of the beam $\beta=\beta (c_1, c_2, \delta)$.  The point $x$ divides the beam $\beta$ into two segments which we denote by $\beta_1$ and $\beta_2$. 

\begin{definition}\label{def-gamma}
By  $b_i=b_i(c_1,c_2,\delta)$, we denote the loop based at the midpoint $x=x(c_1,c_2,\delta)$, obtained by following the half-beam $\beta_i$, then winding once along unit-length horocycle $\partial \Hor_{c_{i}} (1)$, and returning to $x$ back via $\beta_i$. We refer to
$b_1$ and $b_2$ as the lollipops.
\end{definition}
\begin{figure}
\begin{tikzpicture}
\filldraw (2,0) circle (.06); \node at (2.2,-.3) {$c_2$};
\filldraw (-2,0) circle (.06); \node at (-2.2,-.3) {$c_1$};
\draw[thick,dotted,deepblue] (-2,0) circle (1); \draw[thick,dotted,deepred] (2,0) circle (1); \node at (-2,.5) {$\mathcal{H}_{c_1} (1)$}; \node at (2,.5) {$\mathcal{H}_{c_2} (1)$};
\draw[dashed] (-2,0)--(-1,0) (2,0)--(1,0);1
\draw (-1,0)--(1,0);
\draw[deepblue] (-.096,.096)--(0,0)--(-.096,-.096);
\draw[->-,deepblue] (-.096,-.096)--(-.904,-.096); \draw[->-,deepblue] (-.904,.096)--(-.096,.096);
\draw[deepblue] (-.904,.096) arc (5:355:1.1);
\draw[deepred] (.096,.096)--(0,0)--(.096,-.096);
\draw[->-,deepred] (.096,-.096)--(.904,-.096); \draw[->-,deepred] (.904,.096)--(.096,.096);
\draw[deepred] (.904,.096) arc(175:-175:1.1);
\filldraw (0,0) circle (.05) node[above=3pt] {$x$};
\node at (-.5,-.35) {}; \node at (.5,-.35) {}; \node at (-3,1) {$b_1$}; \node at (3,1) {$b_2$}; \node at (-0.5,-.5) {$\beta_1$}; \node at (0.5,-.5) {$\beta_2$};
\end{tikzpicture}
\caption{The lollipops $b_1$ (blue), and $b_2$ (red) represent elements of the fundamental group of $\pi_1(X,x)$}\label{fig-1}
\end{figure}
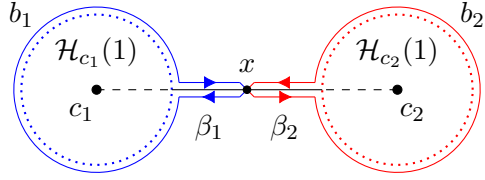

Choose a pointed marking $f:(\Sigma^p_g,\star)\to (X,x)$,  and let $\phi \in \Lab_{g,p,n}$. Then there exist a unique   marked pointed Riemann surface $(X_\phi,x_\phi,f_\phi)$,  and  a holomorphic unbranched covering  $\pi_\phi:X_\phi \to X$, such that the following diagram commutes
\begin{center}
\begin{tikzpicture}[every node/.style={midway}]
  \matrix[column sep={7em,between origins}, row sep={2em}] at (0,0) {
    \node(R) {$(X_\phi,x_\phi)$}  ; & \node(S) {$(\Sigma',\star')$}; \\
    \node(R/I) {$(X,x)$}; & \node (T) {$(\Sigma^p_g,\star)$};\\
  };
  \draw[<-] (R/I) -- (R) node[anchor=east]  {$\pi_\phi$};
  \draw[<-] (R) -- (S) node[anchor=south] {$f_\phi$};
  \draw[->] (S) -- (T) node[anchor=west] {$\pi$};
  \draw[<-] (R/I) -- (T) node[anchor=north] {$f$};
\end{tikzpicture}
\end{center}
where  $\pi:(\Sigma',\star') \to (\Sigma^p_g,\star)$ a pointed cover corresponding to the subgroup $P_n(\phi) \in \Kov_{g,p,n}$. 

Recall that each $\phi \in \Lab_{g,p,n}$ correspond to the equivalence class of a fibre labelled cover  $(\pi,\iota)$ where $\pi:\Sigma'\to \Sigma^p_g$, and $\iota:\pi^{-1}(\star)\to [n]$ is a labelling (that is, we have $\phi=\phi_{[\pi,\iota]}$). Let $(\pi_\phi,\iota_\phi)$ denote the fibre labelled cover where $\pi_\phi:X_\phi \to X$
is the aforementioned holomorphic covering,  and $\iota_\phi:\pi^{-1}_\phi(x)\to [n]$ the labelling defined by  $\iota_\phi=\iota \circ  f^{-1}_\phi$.
We set $x(k)=\iota^{-1}_\phi(k)$. 
\begin{proposition}\label{prop-pair-5}  For each  pair of connected cusps $(c_1,c_2,\delta)$ on $X\in \Mod_{g,p}$,
there exists $a_1,a_2 \in \pi_1(\Sigma^p_g,\star)$ whose powers are mutually distinct, and with the following property. Let  $\phi\in \Lab_{g,p,n}$, and let  $(c_1(k), c_2(k), \delta(k))$ be  the  lift of $(c_1, c_2, \delta)$ (under the covering $\pi_\phi:X_\phi\to X$) whose midpoint is $x(k)$.  Then 
\begin{equation}\label{eq-5}
d_1(k)+d_2(k)=\left|\Orb( \phi(a_1),k)\right| + \left|\Orb( \phi(a_2),k)\right|
\end{equation}
where $d_i(k)$ denotes the branching degree of $\pi_\phi$ at $c_i(k)$.
\end{proposition}
\begin{proof} Note that the branching degree $d_i(k)$ is the smallest positive integer $m$ such that $b_i^{m}$ lifts to a closed loop starting from $x(k)$.  This is equivalent to saying that 
\begin{equation}\label{eq-5-1}
d_i(k)=\left|\Orb( \psi(b_i),k)\right|, \,\,\,\,\,\,\,\,\,\,\, i=1,2,
\end{equation}
where  $\psi:\pi_1(X,x)\to \Sym_n$ is the monodromy homomorphism corresponding to the fibre labelled cover $(\pi_\phi,\iota_\phi)$. Observe that 
$\psi=\phi\circ (f_*)^{-1}$, where $f_*:\pi_1(\Sigma^p_g,\star)\to \pi_1(X,x)$ is the induced isomorphism. Replacing this in (\ref{eq-5-1}) yields the equality 
$$
d_i(k)=\left|\Orb( \phi(a_i),k)\right|
$$
where $a_i=f^{-1}_*(b_i)$. Clearly, $a_1$ and $a_2$ depend only on $(c_1,c_2,\delta)$, and not on $\phi$.
This implies the identity (\ref{eq-5}). The reader can verify that the nontrivial powers of $a_1$ and $a_2$ are
distinct because the parabolic deck transforms of the universal cover induced by $a_1$ and $a_2$ have different  fixed
points.

\end{proof}

\subsection{Proof of Theorem \ref{thm-adam-1}}\label{subsection-proof}

\begin{proposition}\label{prop-pair-4} For each  $X\in \Mod_{g,p}$
there exists a finite collection of pairs  $A\subset  \pi_1(\Sigma^p_g,\star)\times  \pi_1(\Sigma^p_g,\star)$ with the following properties. 
Firstly, if $(a_1,a_2)\in A$ then the powers of $a_1$ and $a_2$ are mutually distinct. Secondly,
suppose   $\phi\in \Lab_{g,p,n}\setminus \Lab^{\Hor}_{n}(X,L)$. Then there are   $(a_1,a_2)\in A$  which are $2\log L$-tangled by $\phi$. 
\end{proposition}
\begin{proof} To each pair of connected cusps $(c_1,c_2,\delta)$ on $X$ we associate the pair $(a_1,a_2)\in\pi_1(\Sigma^p_g,\star)\times  \pi_1(\Sigma^p_g,\star)$ from Proposition \ref{prop-pair-5}.  We define $A$ as the collection of  pairs $(a_1,a_2)$ corresponding to pairs of connected cusps of length at most $2\log L$. Then, the collection $A$ is finite by Proposition \ref{prop-pair-3} (because the corresponding collection of pairs of connected cusps  is finite).

The assumption   $\phi\in \Lab_{g,p,n}\setminus \Lab^{\Hor}_{n}(X,L)$ means that $X_\phi$ does not have the $L$-horoball property. By Proposition \ref{prop-pair-1} there exists a pair of connected cusps $(c'_1,c'_2,\delta')$ on $X_\phi$ whose length is at most $2\log L$. 
Let $(c_1,c_2,\delta)=\pi_\phi(c'_1,c'_2,\delta')$. Then  $(c'_1,c'_2,\delta')=(c_1(k),c_2(k),\delta(k))$ for some $k\in [n]$. 
From Proposition \ref{prop-pair-2} we conclude that 
$$
d_1(k)+d_2(k)+|(c_1,c_2,\delta)| = |(c_1(k),c_2(k),\delta(k))|.
$$
Since  we assume that $|(c_1(k),c_2(k),\delta(k))|\le 2\log L$, it follows that 
\begin{equation}\label{eq-small}
|(c_1,c_2,\delta)|\le 2\log L,
\end{equation}
and
\begin{equation}\label{eq-ville}
d_1(k)+d_2(k)\le 2\log L.
\end{equation}
From Proposition \ref{prop-pair-5} we conclude that for certain $(a_1,a_2)\in A$,  the following holds
$$
d_1(k)+d_2(k)=\left|\Orb( \phi(a_1),k)\right| + \left|\Orb( \phi(a_2),k)\right|.
$$
Combining this with (\ref{eq-ville}) yields  the inequality
$$
\left|\Orb( \phi(a_1),k)\right| + \left|\Orb( \phi(a_2),k)\right|\le 2\log L.
$$
But this means that $a_1$ and $a_2$ are $2\log L$-tangled by $\phi$. This proves the proposition.
\end{proof}

We now complete the proof of Theorem \ref{thm-adam-1}. Fix $X \in \Mod_{g,p}$, and $L>0$. We need to prove the equality (\ref{eq-adam-1}).
From Proposition \ref{prop-pair-4} and Theorem \ref{thm-comb} we conclude that
\begin{align*}
\frac{|\Lab_{g,p,n}\setminus \Lab^{\Hor}_{n}(X,L)|}{|\Lab_{g,p,n}|}&\le \sum_{(a_{1},a_{2})\in A} \frac{\left|\{ \phi\in \Lab_{g,p,n} :\, a_1,a_2 \,\, \text{are ($2\log L$)-tangled by} \,\, \phi \}\right|}{\left|\Lab_{g,p,n}\right|} \\
&\le \frac{C|A|}{n}.
\end{align*}
where $C$ is the constant from Theorem \ref{thm-comb}.
Letting $n\to\infty$ in the previous inequality implies   the equality (\ref{eq-adam-1}).

\section{Random permutations are tangle free}\label{section-comb}
It remains to prove  Theorem \ref{thm-comb-0}. The key statement is that the set of homomorphisms from $\phi\in \Hom_{m,n}$ such that $\phi(w_1)$ and $\phi(w_2)$ have a common fixed point in the set $[n]$ is statistically insignificant compared to the size of the set $\Hom_{m,n}$ (here we must assume that the nontrivial powers of $w_1$ and $w_2$ are distinct.).  We show that each such $\phi$ is carried by an edge labelled graph which  enables us to effectively bound above the number of such homomorphisms $\phi$.

\subsection{Carrier graphs} Given a directed graph $G$ we denote the vertex set of $G$ by $V(G)$, and the set of oriented edges  by $E(G)$. For $e\in E(G)$ we let $\iota(e)$ and $\tau(e)$ denote the initial and terminal vertices of $e$ respectively. We write $\chi(G) = |V(G)| - |E(G)|$ to denote the Euler characteristic of $G$.
Next, we introduce the key definitions of this section.

\begin{definition}\label{definition-carry-0} We say that $(G,h)$ is an edge labelled graph if:
\begin{enumerate}
\item  $G$ is a weakly connected  directed graph,
\vskip .1cm
\item  $h:E(G)\to [m]$ is an edge labelling such that if two edges $e_1,e_2\in E(G)$ have the same initial and terminal vertices  then  $h(e_1)\ne h(e_2)$.
\end{enumerate}
\end{definition}

It turns out that edge labelled graphs are a convenient way of tracking fixed points of permutations.  Let  $\{s_1,\dots, s_m \}$ denote a  generating set of the group $\Free_m$.
\begin{definition}\label{definition-carry} 
We say that  an edge labelled graph $(G,h)$ is  $f$-compatible with  $\phi\in \Hom_{m,n}$ if 
$f :V(G) \to [n]$  is vertex labelling such that 
$$
\phi(s_{h(e)})\big(f(\iota(e))\big)=f(\tau(e))
$$ 
for every  $e\in E(G)$, where $s_{h(e)}$ is the corresponding generator of $\Free_m$.  Furthermore, we say that $(G,h)$ carries $\phi$ if $(G,h)$ is  $f$-compatible with $\phi$ for some  \textbf{injective} vertex labelling $f$.
\end{definition}

The next proposition translates common fixed points of images of $\Free_m$ in the symmetric
group $\Sym_n$  into the language of carrier graphs with negative Euler characteristics.
We postpone its proof until Section \ref{section-comb-1}. 

\begin{proposition}\label{prop-carry}
For every  $w_1, w_2 \in \Free_m$  there exists a constant $C=C(w_1,w_2)$ with the following properties.
Suppose  $\phi \in \Hom_{m,n}$ is such that $\phi(w_1), \phi(w_2)$ have a common fixed point in the set $[n]$.
Then $\phi$ is carried by an edge labelled graph $(G,h)$ which has at most $C$ edges. If in addition we assume that  all non-trivial powers of $w_1$ and $w_2$ are distinct then $\chi(G)<0$.
\end{proposition}

\subsection{The number of carried homomorphisms}

If $w_1$ and $w_2$ are $R$-tangled by some homomorphism $\phi$ then the permutations $\phi(w_1^{R!})$ and $\phi(w_2^{R!})$ have a common fixed point in the set $[n]$. Applying Proposition \ref{prop-carry} to $w_1^{R!},w_2^{R!}\in \Free_m$, 
we find that such $\phi$ is  carried by a suitable edge labelled graph $(G,h)$. 

The next step in the proof of Theorem \ref{thm-comb} is to estimate  the number of homomorphism from $\Hom_{m,n}$ which are carried by a fixed edge  labelled graph.
\begin{lemma}\label{lemma-count}
Let $(G,h)$ be an edge labelled graph. There exists a constant $C=C(G,h,m)$ such that for every integer $n>0$ the following holds
\begin{equation}\label{eq-count}
\frac{|\{ \phi \in\Hom_{m,n}:\,   (G,h) \, \text{carries} \,  \phi  \}|}{|\Hom_{m,n}|}\le Cn^{\chi(G)}
\end{equation}
\end{lemma}
\begin{proof} Fix an edge labelled graph $(G,h)$. The proof of the lemma is based on the following three claims. The first claim is elementary and its proof is left to the reader.

\begin{claim}\label{claim-11} The number of injective vertex labelings $f:V(G) \to [n]$ is equal to 
$$
n(n-1) \cdots (n-|V(G)|+1)=\frac{n!}{(n-|V(G)|)!}.
$$
\end{claim}

\begin{claim}\label{claim-12} Suppose that $(G,h)$ carries some $ \phi \in \Hom_{m,n}$. If $h(e_1)=h(e_2)$  then the implication 
\begin{equation}\label{eq-impl}
\iota(e_1)=\iota(e_2)\,\, \implies \,\, e_1=e_2
\end{equation}
holds. 
\end{claim}
\begin{proof} Since  $(G,h)$ carries  $ \phi$ there exists an injective vertex labelling $f:V(G)\to [n]$ such that   $(G,h)$ is $f$-compatible with  $ \phi$. Suppose $\iota(e_1)=\iota(e_2)$.  It follows from the $f$-compatibility  that  $\phi(s_l)$ sends $f(\iota(e_j))$ to $f(\tau(e_j))$, for  $j=1,2$, where  $l=h(e_1)=h(e_2)$. This implies that $f(\tau(e_1))=f(\tau(e_2))$. Since $f$ is injective we derive the equality $\tau(e_1)=\tau(e_2)$. Thus, we have shown that $e_1$ and $e_2$ have the same  initial and terminal vertices. Combining this with the condition (2) from Definition \ref{definition-carry-0}  shows that $e_1=e_2$, and the implication (\ref{eq-impl}) is proved.
\end{proof}

Next, define $E_l(G)=\{e\in E(G): \, h(e)=l \}$.

\begin{claim}\label{claim-13}  Let  $f:V(G) \to [n]$ be an injective vertex labelling. Then 
\begin{equation}\label{eq-13}
|\{ \phi \in \Hom_{m,n}: \text{$(G,h)$ is $f$-compatible with $\phi$} \}| \le \prod_{l=1}^{m} (n-|E_l(G)|)! 
\end{equation}
\end{claim}
\begin{proof}  We estimate the left hand side in (\ref{eq-13}) as
\begin{equation}\label{eq-13-1}
|\{ \phi \in \Hom_{m,n}: \text{$(G,h)$ is $f$-compatible with $\phi$} \}| \le \prod_{l=1}^{m} X_l
\end{equation}
where $X_l$ is the number of  permutations in $\Sym_n$ which are realised as $\phi(s_l)$ for some 
$\phi$ which is $f$-compatible with $(G,h)$. We now estimate $X_l$.

Define the subset $[n](G,h,f,l)\subset [n]$ by
$$
[n](G,h,f,l)=\{f(\iota(e)) : e \in E_l(G)\}.
$$
If $(G,h)$  is $f$-compatible with some  $\phi$ then the  vertex labelling $f$ specifies the values of the permutation $\phi(s_l)$ on the set $[n](G,h,f,l)$.  Combining  Claim \ref{claim-12} with the assumption that $f$ is injective implies that the set  $[n](G,h,f,l)$ has exactly $|E_l(G)|$ elements. Therefore, there are either zero or
exactly  $(n-|E_l(G)|)!$ possible permutations in $\Sym_n$ that can equal $\phi(s_l)$. That is, 
we established the estimate $X_l \le (n-|E_l(G)|)!$. Replacing this in (\ref{eq-13-1}) proves the claim.
\end{proof}

We are ready to finish the proof of the lemma. The number of homomorphisms $\phi \in \Hom_{m,n}$ which are carried by $(G,h)$ can be estimated above by the product of  two numbers. The first is the number of injective vertex labelings $f$, and the second is the number of $\phi\in \Hom_{m,n}$ which are $f$-compatible with $(G,h)$ for a fixed injective vertex labelling $f$. 
These two numbers we estimated in Claim \ref{claim-11} and Claim \ref{claim-13} respectively, and  we derive the estimate
\begin{equation}\label{eq-susa}
|\{ \phi \in \Hom_{m,n}:\, \text{$(G,h)$ carries $\phi$} \}|\le  \frac{n!}{(n-|V(G)|)!} \prod_{l=1}^{m} (n-|E_l(G)|)! .
\end{equation}
On the other hand,  it is easy to derive (the well known) equality:
\begin{equation}\label{eq-suska}
|\Hom_{m,n}|=(n!)^m.
\end{equation}
Hence, dividing the left-hand side of (\ref{eq-susa}) by the left-hand side of (\ref{eq-suska}) yields the estimate
\begin{align*}
\frac{|\{ \phi \in \Hom_{m,n} :\,  (G,h) \, \text{carries} \,  \phi  \}|}{|\Hom_{m,n}|}&\le\frac{n(n-1) \cdots (n-|V(G)|+1)}{\prod\limits_{l=1}^{m} n(n-1) \cdots (n-|E_l(G)|+1)} \\
&\le C \frac{n^{|V(G)|}}{n^{(|E_1(G)|+\cdots +|E_m(G)|)}}=Cn^{\chi(G)}
\end{align*}
for some $C$ depending only on $G,h$ and $m$. In the last step we used the equality $|E_1(G)|+\cdots+ |E_m(G)|=|E(G)|$. \end{proof}

\subsection{Proof of Theorem \ref{thm-comb}}

It remains to finish the proof of Theorem \ref{thm-comb}. Suppose that $w_1$ and $w_2$ are R-tangled by some $\phi\in \Hom_{m,n}$. 
Then the permutations $\phi(w_1^{R!}), \phi(w_2^{R!})$ have a common fixed point in the set $[n]$. Then Proposition \ref{prop-carry} states that such $\phi$ is carried by an edge labelled graph $(G,h)$ with at most $C_1$ edges, and with the negative Euler characteristic. Here $C_1$ is the constant from Proposition \ref{prop-carry} which depends only on $w_1^{R!}$ and $w_2^{R!}$.
This implies the following estimate
\begin{equation}\label{eq-AB}
\left|\{ \phi\in \Hom_{m,n} : w_1,w_2 \,\, \text{are R-tangled by} \,\, \phi \}\right|\le A  B
\end{equation}
where $A$ is the number of homomorphism $\phi$ carried by a fixed edge labelled graph $(G,h)$ with $\chi(G)<0$, and $|E(G)|\le C_1$,
and $B$ is the number of edge labelled graphs $(G,h)$ with $\chi(G)<0$, and $|E(G)|\le C_1$. 

The number of graphs with at most $C_1$ edges is a finite number depending only on $C_1$. For each such graph $G$, the number of edge labelings $h:E(G)\to [m]$ depends only on $|E(G)|$ and $m$. Thus,  the number $B$ depends only $C_1$ and $m$. We conclude that $B$ depends only on $w_1,w_2,R,m$.

On the other hand, for a fixed $(G,h)$ we estimate $A$ using  Lemma \ref{lemma-count} :
\begin{align*}
\frac{A}{\left|\Hom_{m,n}\right|}&= \frac{|\{ \phi \in\Hom_{m,n} :\,   (G,h) \, \text{carries} \,  \phi  \}|}{|\Hom_{m,n}|}\\
&\le C_2n^{\chi(G)}\le \frac{C_2}{n},
\end{align*}
where $C_2$  is the constant from Lemma \ref{lemma-count} depending only on $G,h$ and $m$. In the last step we used that $\chi(G)<0$. Returning this to (\ref{eq-AB}) yields
$$
\frac{\left|\{ \phi\in \Hom_{m,n} : w_1,w_2 \,\, \text{are R-tangled by} \,\, \phi \}\right|}{\left|\Hom_{m,n}\right|}\le \frac{C}{n}
$$
for $C=BC_2$. We have proved the theorem.

\section{Proof of Proposition \ref{prop-carry}}\label{section-comb-1}

The first step is to construct an edge labelled graph $(G,h)$ which is $f$-compatible with $\phi$ without the requirement that $f$ is injective. This is the content of Lemma \ref{lemma-g-1} below. In the endgame we modify $(G,h)$ to a new edge labelled graph $(G_1,h_1)$ which carries $\phi$.

Let $Y$ be the directed graph which has a single vertex $y$, and $m$ oriented loops (petals) which we denote by 
$\alpha_1,\dots,\alpha_m$. There is an obvious isomorphism $\pi_1(Y,y)\to \Free_m$ defined by $\alpha_j\to s_j$, $j\in [m]$.

\begin{figure}
\begin{align*}
\phi \colon \begin{cases} s_1 &\mapsto (15)(687)\\ s_2 &\mapsto (172569)\\ s_3 &\mapsto (1934)(58) \end{cases}\\
\end{align*}
$$
w_1= s_1^{-1}s_2^{-1}s_1s_2, \,\,\,\,\,\,\,\,\,\,\,\,\,\,\,\,  \phi(w_1)= (29)(5678)  
$$
$$
w_2 = s_1^{-1}s_3^{-1}s_1s_2^{-1}s_3,  \,\,\,\,\,\,\,\,\,\,\,  \phi(w_2)= (276458)
$$
\begin{tikzpicture}[scale=2] 
\draw[->>-,deepgreen] (0,0)--(-1,1);\node[rotate=-45] at (-.65,.35) { \footnotesize{$h(e_{1,4})=2$}};
\draw[->-,deepblue] (-1,1)--(-2,0);\node[rotate=45] at (-1.7,.6) { \footnotesize{$h(e_{1,3})=1$}};
\draw[->>-,deepgreen] (-1,-1)--(-2,0);\node[rotate=-45] at (-1.7,-.6) { \footnotesize{$h(e_{1,2})=2$}};
\draw[->-,deepblue] (0,0)--(-1,-1);\node[rotate=45] at (-.65,-.35) { \footnotesize{$h(e_{1,1})=1$}};
\draw[->>>-,deepred] (0,0)--(.691,.951);\node[rotate=54] at (.25,.6) { \footnotesize{$h(e_{2,5})=3$}};
\draw[->>-,deepgreen] (1.809,.588)--(.691,.951);\node[rotate=-18] at (1.35,.95) { \footnotesize{$h(e_{2,4})=2$}};
\draw[->-,deepblue] (1.809,.588)--(1.809,-.588);\node[rotate=-90] at (2,0) { \footnotesize{$h(e_{2,3})=1$}};
\draw[->>>-,deepred] (.691,-.951)--(1.809,-.588);\node[rotate=18] at (1.35,-.95) { \footnotesize{$h(e_{2,2})=3$}};
\draw[->-,deepblue] (0,0)--(.691,-.951);\node[rotate=-54] at (.25,-.6) { \footnotesize{$h(e_{2,1})=1$}};
\node[circle,draw=black,fill=white] at (0,0) {$1$};\node at (.3,0) {$v_0$};
\node[circle,draw=black,fill=white] at (-1,1) {$7$};\node[circle,draw=black,fill=white] at (-2,0) {$6$};\node[circle,draw=black,fill=white] at (-1,-1) {$5$};
\node[circle,draw=black,fill=white] at (.691,.951) {$9$};\node[circle,draw=black,fill=white] at (1.809,.588) {$6$};\node[circle,draw=black,fill=white] at (1.809,-.588) {$8$};\node[circle,draw=black,fill=white] at (.691,-.951) {$5$};
\end{tikzpicture}
\caption{Both permutations $\phi(w_1), \phi(w_2)$ fix $1$. Vertices represent the orbits of the common fixed point 1 under  $\phi(w_1)$ and $\phi(w_2)$, and are labelled according the vertex labelling $f$.}\label{fig-2}
\end{figure}
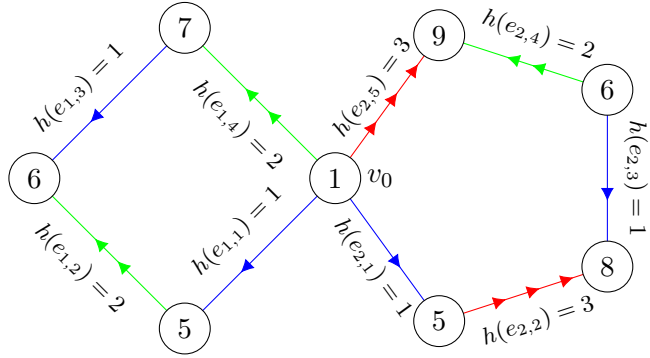

\begin{lemma}\label{lemma-g-1} Let $w_1, w_2 \in \Free_m$, and suppose  $\phi \in \Hom_{m,n}$ is such that the permutations  
$\phi(w_1), \phi(w_2)$ have a common fixed point in the set $[n]$.
Then there exists an edge labelled graph $(G,h)$ with the following properties:
\begin{enumerate}
\item the number of edges $|E(G)|$ is equal to the sum of the (reduced) word lengths of $w_1$ and $w_2$ in $\Free_m$,  
\vskip .1cm
\item there exists a (non necessarily injective) vertex labelling $f:V(G)\to [n]$ such that $(G,h)$ is  $f$-compatible with $\phi$,
\vskip .1cm
\item there exists a graph morphism $\mu:G\to Y$ such that (for a suitable $v_0 \in V(G)$) the group $\mu_*(\pi_1(G,v_0))$ is equal to the subgroup of $\pi_1(Y,y)$ generated by $w_1$ and $w_2$,
\vskip .1cm
\item the equivalence  
$$
h(e_1)=h(e_2) \,\, \iff \,\, \mu(e_1)=\mu(e_2)
$$ 
holds for every $e_1,e_2\in E(G)$. 
\end{enumerate}
\end{lemma}
\begin{remark} The significance of the condition $(3)$ is that implies  $\chi(G)<0$ provided  the group generated by 
$w_1$ and $w_2$ is not Abelian. 
\end{remark}
\begin{remark} Figure \ref{fig-2} illustrates the construction of the edge labelled graph $(G,h)$ and the corresponding vertex labelling.
\end{remark}

\begin{proof} Recall that $\{s_1,\dots,s_m\}$ is the generating set for $\Free_m$. Therefore, each $w \in \Free_m$ can be written as a (reduced) word in the corresponding alphabet. In particular, the elements $w_1$ and $w_2$ are spelled as:
\begin{equation}\label{eq-word}
w_i=s^{\sigma_{i,1}}_{l_{i,1}} \,  s^{\sigma_{i,2}}_{l_{i,2}} \cdots s^{\sigma_{i,k_{i}}}_{l_{i,k_{i}}}
\end{equation}
where $\sigma_{i,t}\in \{1,-1\}$.

Let $Z$ be a directed graph which has a single vertex $z$, and $2$ oriented loops (petals) which we denote by 
$\beta_1$ and $\beta_2$.  The petal $\beta_i$ is subdivided into  oriented edges $e_{i,t}$, where  $1\le t \le k_i$, so that 
$e_{i,t}$ has the same orientation as $\beta_i$ if and only of $\sigma_{i,t}=1$. The resulting directed graph is denoted by $G$. 
Then $|E(G)|$  equals  the sum of word lengths of $w_1$ and $w_2$ which  confirms  (1).

Next, we define the morphism $\mu:G\to Y$ by letting 
\begin{equation}\label{eq-word-1}
\mu(e_{i,t})=\alpha_{l_{i,t}}. 
\end{equation}
Let  $v_0\in V(G)$ be the vertex arising from the vertex $z$ of the rose $Z$.  Then  $\pi_1(G,v_0)$ is a  free group of rank two  generated by the loops $\beta_1$ and $\beta_2$. From  (\ref{eq-word}) we conclude  $\mu_*(\beta_i)=w_i$, where $\mu_*:\pi_1(G,v_0)\to \Free_m$ is the induced homomorphism.  This confirms   (3).

Define the edge labelling $h:E(G)\to [m]$ by letting 
\begin{equation}\label{eq-word-2}
h(e_{i,t})=l_{i,t}.  
\end{equation}
Any two vertices in $V(G)$ are connected by at most one edge. Therefore, $(G,h)$ is an edge labelled graph in the sense of Definition \ref{definition-carry-0}. Moreover, combining definitions (\ref{eq-word-1}) and (\ref{eq-word-2}) proves (4).

It remains to construct a vertex labelling $f:V(G)\to [n]$ such that $(G,h)$ is  $f$-compatible with $\phi$.  Denote by $k\in [n]$ the common fixed point of the permutations by $\phi(w_1)$ and  $\phi(w_2)$. Set $f(v_0) = k$, and propagate
the vertex labelling $f$ to all other vertices so as to satisfy the condition 
$$
\phi(s_l) \big(f( \iota(e) ) \big)=f(\tau(e))
$$
for each $e$, where $h(e)=l$. Since $k$ is fixed by both permutations  $\phi(w_1)$ and  $\phi(w_2)$, such vertex labelling $f$ is well defined and  is in fact unique. Thus, we have shown that $(G,h)$ is $f$-compatible with $\phi$ which confirms (2). 
\end{proof}

\subsection{The endgame}

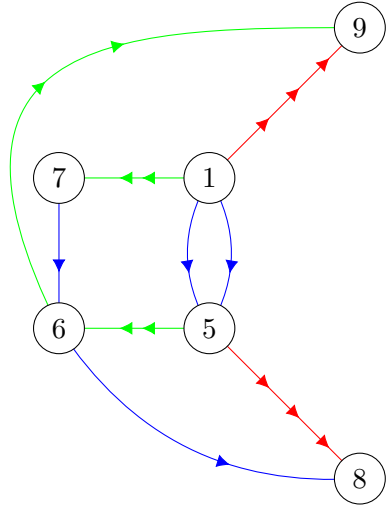
\begin{figure}
\begin{tikzpicture}[scale=2]
\draw[->>-,deepgreen] (0,0)--(-1,0);
\draw[->-,deepblue] (-1,0)--(-1,-1);
\draw[->>-,deepgreen] (0,-1)--(-1,-1);
\draw[->-,deepblue] (0,0) to[bend right] (0,-1);
\draw[->-,deepblue] (0,0) to[bend left] (0,-1);
\draw[->>>-,deepred] (0,0)--(1,1);
\draw[->>>-,deepred] (0,-1)--(1,-2);
\draw[->-,deepblue] (-1,-1) to[bend right] (1,-2);
\draw[->>-,deepgreen] (-1,-1) ..controls (-2,1) and (-.5,1).. (1,1);
\node[circle,draw=black,fill=white] at (0,0) {$1$};
\node[circle,draw=black,fill=white] at (-1,0) {$7$};
\node[circle,draw=black,fill=white] at (-1,-1) {$6$};
\node[circle,draw=black,fill=white] at (0,-1) {$5$};
\node[circle,draw=black,fill=white] at (1,1) {$9$};
\node[circle,draw=black,fill=white] at (1,-2) {$8$};
\end{tikzpicture}
\caption{The graph $G'$ from Section 7.1. The pairs of vertices with the labels 5 and 6 were identified by the relation $\sim_V$. The graph $G_1$ is then obtained from $G'$ by replacing the two blue edges of label 1 between the vertices 1 and 5 with a single edge. Then $\chi(G_1)=-2$.}\label{fig-3}
\end{figure}

Given $\phi\in \Hom_{m,n}$, in Lemma \ref{lemma-g-1} we have constructed an edge labelled graph $(G,h)$ which is $f$-compatible with $\phi$, where $f:V(G)\to [n]$ is a labelling. We define a new graph $G_1$, together with a graph morphism $\rho:G\to G_1$, as follows (see Figure \ref{fig-3}).

Define the equivalence relation $\sim_V$ on the vertex set $V(G)$ by letting $v_1\sim v_2$ if and only if $f(v_1)=f(v_2)$. Denote by $G'$ the corresponding quotient graphs such that $V(G')=V(G)/\sim$. Then, the graph $G_1$ is obtained from $G'$ by a sequence of  foldings where at each stage we identify two edges $e_1$ and $e_2$ if they have the same initial and terminal vertices, and the same label (that is, $h(e_1)=h(e_2)$). By $\rho:G\to G_1$ we denote the resulting morphism. 

Given  $e_1\in E(G_1)$, there exists $l\in [m]$ such that $h(e)=l$ for every $e\in \rho^{-1}(e_1)\subset E(G)$. We define the new edge labelling $h_1:E(G_1)\to [m]$ by letting $h_1(e_1)=l$. By construction $(G_1,h_1)$ is an edge labelled graph  in the sense of Definition \ref{definition-carry-0}. 

Likewise, for each  $v_1 \in V(G_1)$  there exists $k\in [n]$ such that $f(v)=k$ for every $v\in \rho^{-1}(v_1)\subset V(G)$. This enables us to  define the vertex labelling $f_1:V(G_1)\to [n]$ by letting $f_1(v_1)=k$. By construction $f_1$ is injective (since the vertex $v_1\in V(G_1)$ is  the equivalence class consisting of all vertices from $V(G)$ that are mapped to $k$). Moreover, since $(G,h)$ is $f$-compatible with $\phi$ we conclude that $(G_1,h_1)$ is $f_1$-compatible with $\phi$.

Thus, we have constructed an edge labelled graph $(G_1,h_1)$ which carries $\phi$. It remains to compute its Euler characteristic. It follows from the condition (4) in Lemma \ref{lemma-g-1} that the morphism $\mu:G\to Y$ factors through $G_1$. That is, there exist a morphism $\mu_1:G_1\to Y$ so that $\mu=\mu_1 \circ \rho$. Then the condition (3) from Lemma  \ref{lemma-g-1} implies the inclusion 
\begin{equation}\label{eq-ram}
W= \mu_*\big(\pi_1(G,v_0)\big) \le (\mu_1)_*\big(\pi_1(G_1,v_1)\big)
\end{equation}
where $W$ is the subgroup of $\pi_1(Y,y)$ generated by the words $w_1$ and $w_2$ (here $v_1=\rho(v_0)$). If $w_1$ and $w_2$ have no non-trivial  powers in common then $W$ is a free group of rank $2$. Combining this together with (\ref{eq-ram}) yields the inequality $\chi(G)<0$. 

Therefore, we have proved that the labelled graph $(G_1,h_1)$ carries $\phi$. On the other hand, the obvious inequality $|E(G_1)|\le |E(G)|$, and the condition (1) from  Lemma  \ref{lemma-g-1}, complete the proof of Proposition \ref{prop-carry}.

\end{document}